\documentclass{amsart}

\usepackage{amsmath} 
\usepackage{amssymb}
\usepackage{mathrsfs}

\newtheorem{theorem}{Theorem}[section] 
\newtheorem{claim}[theorem]{Claim}

\newtheorem{conclusion}[theorem]{Conclusion}

\theoremstyle{definition}
\newtheorem{definition}[theorem]{Definition}

\newtheorem{observation}[theorem]{Observation}

\theoremstyle{remark}
\newtheorem{question}[theorem]{Question}
\newtheorem{remark}[theorem]{Remark}
\newtheorem{notation}[theorem]{Notation}
\newtheorem{discussion}[theorem]{Discussion}
\newtheorem{context}[theorem]{Context}

\newcommand{\rest}{{\restriction}}

\newcommand{\wilog}{{\rm without loss of generality}}

\newcommand{\then}{{\underline{then}}}
\newcommand{\when}{{\underline{when}}}
\newcommand{\Then}{{\underline{Then}}}

\newcommand{\mn}{{\medskip\noindent}}
\newcommand{\sn}{{\smallskip\noindent}}

\newcommand{\bbB}{{\mathbb B}}

\newcommand{\cC}{{\mathscr C}}

\newcommand{\bbD}{{\mathbb D}}
\newcommand{\cD}{{\mathscr D}}

\newcommand{\cE}{{\mathscr E}}

\newcommand{\cF}{{\mathscr F}}

\newcommand{\cI}{{\mathscr I}}

\newcommand{\bbL}{{\mathbb L}}

\newcommand{\gC}{{\mathfrak C}}

\newcommand{\cP}{{\mathscr P}}

\newcommand{\cT}{{\mathscr T}}
 
\newcommand{\cU}{{\mathscr U}}

\newcommand{\cf}{{\rm cf}}

\newcount\skewfactor
\def\mathunderaccent#1#2 {\let\theaccent#1\skewfactor#2
\mathpalette\putaccentunder}
\def\putaccentunder#1#2{\oalign{$#1#2$\crcr\hidewidth
\vbox to.2ex{\hbox{$#1\skew\skewfactor\theaccent{}$}\vss}\hidewidth}}

\newenvironment{PROOF}[2][\proofname.]
   {\begin{proof}[#1]\renewcommand{\qedsymbol}{\textsquare$_{\rm #2}$}}
   {\end{proof}}

\begin{document}

\title {No limit model in inaccessible}
\author {Saharon Shelah}
\address{Einstein Institute of Mathematics\\
Edmond J. Safra Campus, Givat Ram\\
The Hebrew University of Jerusalem\\
Jerusalem, 91904, Israel\\
 and \\
 Department of Mathematics\\
 Hill Center - Busch Campus \\ 
 Rutgers, The State University of New Jersey \\
 110 Frelinghuysen Road \\
 Piscataway, NJ 08854-8019 USA}
\email{shelah@math.huji.ac.il}
\urladdr{http://shelah.logic.at}
\thanks{The author would like to thank the Israel Science Foundation for
partial support of this research (Grant No.242/03). Publication 906.\\ 
The author thanks Alice Leonhardt for the beautiful typing.}

\date{January 4, 2011}

\dedicatory {Dedicated to Michael Makkai}

\begin{abstract}
Our aim is to improve the negative results i.e., 
non-existence of limit models, and the failure of the generic
pair property from \cite{Sh:877} to inaccessible $\lambda$ as promised
there.  In \cite{Sh:877}, the negative results were obtained only for 
non-strong limit cardinals.
\end{abstract}

\maketitle
\numberwithin{equation}{section}
\setcounter{section}{-1}

\section {Introduction} 

Let $\lambda = \lambda^{<\lambda} > \kappa$ be regular cardinals.  A
complete first order theory $T$ may have (some variant of)
$(\lambda,\kappa)$-limit model, which, if exists, is unique, see
history in \cite{Sh:877} and Definition \ref{0.26}.  
There we prove existence for the theory of
linear order and non-existence for first order theories which are
strongly independent and then just independent and even the parallel for
$\kappa=2$ (one direction of the so-called genreic pair conjecture).  Those
non-existence results in \cite{Sh:877} 
were for $\lambda = 2^\kappa$, here we deal with strongly
inaccessible $\lambda$.  In \cite{Sh:900} there are existence results
but for $\lambda$ measurable, and we promise there the non-existence
results for $\lambda$ strongly inaccessible as complimentary results. 

Let $\lambda$ be strongly inaccessible $(> |T|)$ such that $\lambda^+ =
2^\lambda$; this for transparency only.

Here in \S1 we prove that for strongly independent $T$ (see Definition
\ref{0.7}), a strong version of the generic pair conjecture (see Definition
\ref{0.14}(2)) holds.    
We also prove the non-existence of $(\lambda,\kappa)$-limit
models, a related property (for all versions of ``limit model").

In \S2, we also prove this even for independent $T$.  The use of
$\lambda^+ = 2^\lambda$ is just to have a more transparent formulation
of the conjecture.
See more on the generic pair conjecture for dependent $T$ in
\cite{Sh:950}.

We thank Itay Kaplan for much helpful criticism.
\bigskip

\noindent
\begin{notation}
\label{0.4} 
1) ${\cD}_\lambda$ is the club filter on
$\lambda$ for $\lambda$ regular uncountable.

\noindent
2) $S^\lambda_\kappa = \{\delta < \lambda$\,:\,cf$(\delta)=\kappa\}$.

\noindent
3) For a limit ordinal $\delta$ let ${\cP}^{\text{ub}}(\delta) =
\{{\cU}:{\cU}$ is an unbounded subset of $\delta\}$. 

\noindent
4) $T$ denotes a complete first order theory.

\noindent
5) For a model $M,\varphi(\bar x,\bar y) \in \bbL(\tau_M)$ and $\bar d
 \in {}^{\ell g(\bar y)}M$, let $\varphi(M,\bar d) = \{\bar c \in
 {}^{\ell g(\bar x)}M:M \models \varphi[\bar c,\bar d]\}$.

\noindent
6) $\bold S^n(A,M) = \{\text{tp}(\bar b,A,N):M \prec N$ and $\bar b
 \in {}^n N\}$ where tp$(\bar c,A,N) = \{\varphi(\bar x,\bar
 a):\varphi(\bar x,\bar y) \in \bbL(\tau_M),\bar a
 \in {}^{\ell g(\bar y)}A$ and $M \models \varphi[\bar c,\bar a]\}$.

\noindent
7) $\bold S^n(M) = \bold S^n(M,M)$ and $\bold S^{< \omega}(M) 
= \cup\{\bold S^n(M):n \in \omega\}$.

\noindent
8) $\beth_\alpha(\lambda) = \lambda + \Sigma\{2^{\beth_\beta(\lambda)}:\beta
< \alpha\}$ and $\beth_\alpha = \beth_\alpha(\aleph_0)$.
\end{notation}

Recall (as in \cite[2.3]{Sh:877})
\begin{definition}
\label{0.7}  
1) $T$ has the strong independence
property (or is strongly independent) \when \,:  some 
$\varphi(\bar x,\bar y) \in \bbL(\tau_T)$ has it, where:

\noindent
2)  $\varphi(\bar x,\bar y) \in \bbL(\tau_T)$ has the strong 
independence property for $T$ \when \, for every $n < \omega$, model 
$M$ of $T$ and pairwise disjoint finite $\bold I_1,\bold I_2
\subseteq {}^{\ell g(\bar y)}(M)$ 
for some $\bar a \in {}^{\ell g(\bar x)}M$ we have $\ell \in
\{1,2\} \wedge \bar b \in \bold I_\ell \Rightarrow M \models \varphi[\bar a,
\bar b_\ell]^{\text{if}(\ell =2)}$. 
\end{definition}

\begin{remark}
1) Elsewhere we use $\varphi(x,y)$, i.e. the $x$ and $y$ are
singletons, but the proofs are not affected.

\noindent
2) Also we may restrict ourselves to $\bold I_1,\bold I_2 \subseteq \psi(M,\bar
d)$ where $\psi \in \bbL(\tau_T)$ such that $\psi(M,\bar d)$ is
infinite, and we may restrict ourselves to $\bold I_1,\bold I_2$ such
that every $\bar b \in \bold I_1 \cup \bold I_1$ realizing a fixed 
non-algebraic type $p
\in \bold S^m(A,M)$ with $M$ being $(|A|^+ + \aleph_0)$-saturated.
The results are not really affected. 
\end{remark}

\begin{question}
\label{0q.37}  
1) Assume $\lambda_2 = 
\lambda^{< \lambda_1}_2 \ge \lambda_1 > |T|,T$ a complete first order
dependent theory.
Is the theory $T^*_{\lambda_1,\lambda_2}$ a dependent theory or at
least when is $T^*_{\lambda_1,\lambda_2}$ a dependent theory? where
\mn
\begin{enumerate}
\item[$(a)$]   $T^*_{\lambda_1,\lambda_2} = 
\text{ Th}(K^+_{\lambda_1,\lambda_2})$ where
\sn
\item[$(b)$]  $K^+_{\lambda_1,\lambda_2} = \{(N,M):M \text{ is a }
\lambda_1$-saturated model of $T$ of cardinality $\lambda_2,N$  a 
$\lambda^+_2$-saturated elementary extension of $M\}$.
\end{enumerate}
\mn
2) Similarly for other properties of $T^*_{\lambda_1,\lambda_2}$;
note\footnote{let $(N_\ell,M_\ell) \in K^+_{\lambda_1,\lambda_2}$ for
$\ell =1,2$ and let $f_*$ be an isomorphism from $M_1$ onto $M_2$ and
let $\cF = \{f:f$ is a $(N_1,N_2)$-elementary mapping extending $f_*$
of cardinality $\le \lambda_1\}$.  Now we can prove that any $f
\in \cF$ preserve satisfaction for first order formulas.}
that this theory is complete if $\lambda_1 = \lambda_2$.

\noindent
3) When can we prove that $T^*_{\lambda_1,\lambda_2}$ does not 
depend on the cardinals at least for many pairs?
\end{question}

\begin{remark}
1) Concerning failure of \ref{0q.37}(1) see Kaplan-Shelah
   \cite{KpSh:946}.

\noindent
2) Any solution of the generic pair conjecture answers positively 
\ref{0q.37}(3) for dependent $T$ in the relevant cases. 

\noindent
3) It is known that in \ref{0q.37}(1) if $T$ extends PA or ZFC then
in $T^* = \text{ Th}(N,M)$ we can interpret the second order theory
of $\lambda_2$.
\end{remark}

But may well be that as in Baldwin-Shelah \cite{BlSh:156}
\begin{question}
\label{0q.44}  
Assume $|T| < \kappa \le \lambda_1 \le
\lambda_2 = \lambda^{< \lambda_1}_2,T$ a complete first order theory.  For
which $T$'s can we interpret in $M \in K^+_{\lambda_1,\lambda_2}$
 a model of PA of cardinality $\ge \lambda_1$ by first order formula
\underline{or} just an
$\bbL_{\infty,\kappa}(\tau_T)$-formulas with
parameters, the intention is that we assume $\lambda_2$ is enough
larger than $\lambda_1$
which is large enough than $|T|$; if $2^\kappa \ge \lambda_1$ this is trivial.
\end{question}

Recall (from \cite[0.2]{Sh:877})
\begin{definition}
\label{0.14}   0) Let EC$_\lambda(T)$ be the class of models 
$M$ of (the first order) $T$ of cardinality $\lambda$.  Let
EC$_{\lambda,\kappa}(T)$ be the class of $\kappa$-saturated models
$M \in \text{ EC}_\lambda(T)$.

\noindent
1) Assume $\lambda > |T|$, (we usually assume $\lambda =
\lambda^{< \lambda}$) and $2^\lambda = \lambda^+,M_\alpha \in \text{
EC}_\lambda(T)$ is $\prec$-increasing continuous for $\alpha < 
\lambda^+$ with $M= \cup\{M_\alpha:\alpha < \lambda^+\} \in 
\text{ EC}_{\lambda^+}(T)$, and $M$ is saturated. 
The generic pair property (for $T,\lambda$)
says that for some club
$E$ of $\lambda^+$ for all pairs $\alpha < \beta$ of ordinals from $E$
of cofinality $\lambda,(M_\beta,M_\alpha)$ has the same isomorphism
type (we denote this property of $T$ by Pr$^2_{\lambda,\lambda}(T)$).

\noindent
2) The generic pair conjecture for $\lambda > \aleph_0$, (usually
$\lambda = \lambda^{< \lambda}$) 
such that $2^\lambda = \lambda^+$ says that for any complete
first order $T$ of cardinality $< \lambda,T$ is
dependent iff it has the generic pair property for $\lambda$.

\noindent
3) Let $\bold n_{\lambda,\kappa}(T)$ be 
min$\{|\{M_\delta/\cong \,:\,\delta \in E$ has
cofinality $\kappa\}|:E$ a club of $\lambda^+\}$ for $\lambda$ and
$\bar M = \langle M_\alpha:\alpha < \lambda^+\rangle$ as above and
$\kappa = \text{ cf}(\kappa) \le \lambda$; clearly the choice of
$\bar M$ is immaterial.
\end{definition}

\begin{remark}
\label{0.22} 
1) Note that to say $\bold n_{\lambda,\kappa}(T)=1$ is a way to
say that $T$ has (some variant of) a $(\lambda,\kappa)$-limit model,
see \ref{0.26} below.  There are other variants of the Definition of limit.

\noindent
2) Recall that we conjecture that for $\lambda = \lambda^{< \lambda} >
\kappa = \text{ cf}(\kappa) > |T|,2^\lambda = \lambda^+$ we have
$\bold n_{\lambda,\kappa}(T) =1 \Leftrightarrow \bold n_{\lambda,\kappa}(T)
< 2^\lambda \Leftrightarrow T$ is dependent.  The use of $``\lambda^+ =
2^\lambda"$ is just for clarity.  See more in \cite{Sh:877},
\cite{Sh:900}, \cite{Sh:950}.

\noindent
3) Recall that if $\lambda = \kappa = \lambda^{<\lambda}$, then $\bold
   n_{\lambda,\lambda}(T)=1$ means $T$ has a unique saturated
   model; (and parallely if $\lambda > \cf(\lambda) = \kappa,\lambda$
   strong limit).  So we concentrate on the case 
$\bold n_{\lambda,\kappa}(T)=1$ where $\kappa < \lambda$.
\end{remark}

\begin{definition}
\label{0.26}
We define when $M_*$ is a $(\lambda,\kappa)$-limit model of $T$
   where $\lambda \ge \kappa = \cf(\kappa)$ and $\lambda \ge |T|$
\mn
\begin{enumerate}
\item[$(A)$]  if $2^\lambda = \lambda^+$ this means $\bold
n_{\lambda,\kappa}(T)=1$
\sn
\item[$(B)$]  in general it means that: letting $K_\lambda = \{M:M$ is
a model of $T$ with universe an ordinal $\in [\lambda,\lambda^+)\}$,
for some function $\bold F$ with domain $K$ and satisfying $M \prec
\bold F(M) \in K$ we have:
\sn
\begin{enumerate}
\item[$\oplus$]  if $M_\alpha \in K$ for $\alpha < \lambda^+$ is
$\prec$-increasing continuous and $\alpha < \lambda \Rightarrow \bold
F(M_{\alpha +1}) \prec M_{\alpha + 2}$ then for some club $E$ of
$\lambda^+$ we have:

$\delta \in E \wedge \cf(\delta) = \kappa \Rightarrow M_\delta \cong
M_*$.
\end{enumerate}
\end{enumerate}
\end{definition}

\section {Strongly independent $T$} 

\begin{context}  1) $T$ is a fixed first order complete theory and ${\gC} 
= {\gC}_T$ a monster for it; for notational simplicity $\tau_T$ is relational.

\noindent
2) We let $\lambda$ be a regular uncountable cardinal $> |T|$; we deal
   mainly with strongly inaccessible $\lambda$.
\end{context}

Here for $\lambda$ strongly inaccessible and (complete first order)
$T$ with the strong independence property (of cardinality $< \lambda$) we
prove the non-existence of $(\lambda,\kappa)$-limit models for $\kappa
= \text{ cf}(\kappa) < \lambda$ (in Theorem \ref{inv.14}) and the
generic pair conjecture for $\lambda$ and $T$, in Theorem
\ref{inv.24} (which shows non-isomorphism).  Recall that the generic
pair property speaks on the isomorphism type of pairs of models.

Definition \ref{inv.7} gives us a more constructive invariant of
$(N,M)/\cong$.  Unfortunately it seemed opaque how to manipulate it so
we shall use a related but different version, the one from Definition
\ref{inv.10}.  Naturally it concentrates on types in one formula
$\varphi(\bar y,\bar x)$ witnesssing the strong independence property.  But
mainly gives the pair $(N,M)$ an invariant $\langle {\cP}_\delta:
\delta < \lambda\rangle/{\cD}_\lambda$ where ${\cP}_\delta 
\subseteq {\cP}({\cP}(\delta))$.  Now always 
$|{\cP}_\delta| \le 2^{|\delta|}$ and it is easily computable from one
${\cP} \subseteq {\cP}(\delta)$, in fact from the invariant
inv$_4(M,N)$ from Definition \ref{inv.7}, but in our proofs its use
is more transparent.  It has monotonicity property and we can increase
it.

We need different but similar version for the proof of non-existence
of $(\lambda,\kappa)$-limit models.

\begin{definition}
\label{inv.7}  
1) Let ${\cE}^*_T$ be the
following two-place relation on $\{(M,\bold P):M \models T$ and $\bold P
\subseteq \bold S^{< \omega}(M)\}$; let $(M_1,{\bold P}_1) 
{\cE}^*_T(M_2,\bold P_2)$ iff there is an isormorphism $h$ from $M_1$
onto $M_2$ mapping $\bold P_1$ onto $\bold P_2$. 

\noindent
2) For models $M \subseteq N$ we define (the important case is $M
\prec N \models T$)
\mn
\begin{enumerate}
\item[$(a)$]  $\text{inv}_1(M,N) = 
\{p \in \bold S^{< \omega}(M):p \text{ is realized in } N\}$
\sn
\item[$(b)$]  $\text{inv}_2(M,N) = (M,\text{inv}_1(M,N))/{\cE}^*_T$.
\end{enumerate}
\mn
3) If $M \prec N$ are models of $T$ such that the universe of $N$ is 
$\subseteq \lambda$, recalling ${\cD}_\lambda$ is the club
filter on $\lambda$, let:
\mn
\begin{enumerate}
\item[$(a)$]   for any ordinal $\delta < \lambda$

inv$_3(\delta,M,N) = (M \restriction \delta,
\{p \in \bold S^{< \omega}(M \restriction \delta):p$ is 
realized by some 

\hskip75pt sequence from $N \rest \delta\})/{\cE}^*_T)$
\sn
\item[$(b)$]   $\text{inv}_4(M,N) =
\langle\text{inv}_3(\delta,M,N):\delta < \lambda\rangle/{\cD}_\lambda$.
\end{enumerate}
\mn
4) If $M \prec N$ are models of $T$ of cardinality $\lambda$ then
   inv$_4(M,N)$ is inv$_4(f(M),f(N))$ for every one-to-one function
   $f$ from $N$ into $\lambda$ (equivalently some $f$, see
\ref{inv.8}(1),(2) below)
\end{definition}

\begin{observation}
\label{inv.8}  
0) In Definition \ref{inv.7}(3) for a club of $\delta$'s below $\lambda$ we
   have $M \rest \delta \prec M$ and $N \rest \delta \prec N$ and so
   $M \rest \delta \prec N \rest \delta \models T$.

\noindent
1) Concerning Definition \ref{inv.7}(3), if 
$M \prec N$ are models of $T$ of cardinality
$\lambda$ and $f_1,f_2$ are one-to-one functions from $N$ into
$\lambda$ \then \, inv$_4(f_1(M),f_1(N)) = \text{
inv}_4(f_2(M),f_2(N))$ using the definition \ref{inv.7}(3)(b).

\noindent
2) Definitions \ref{inv.7}(3), \ref{inv.7}(4) are compatible and in
\ref{inv.7}(4), ``some $f$ such that $f$ is a one-to-one function from
$N$ to $\lambda$" is equivalent to ``every $f$ such that..." 
\end{observation}

\begin{PROOF}{\ref{inv.8}}
Straight, e.g. (this argument will be used several times).

\noindent
1) Let $E = \{\delta < \lambda:\delta$ is a limit ordinal such that $M
   \rest \delta \prec M,N \rest \delta \prec N$ and Rang$(f_\ell \rest
   \delta) = \text{ Rang}(f_\ell) \cap \delta$ for $\ell=1,2\}$.  So
   $E$ is a club of $\lambda$ and $\delta \in E \Rightarrow f_2 \circ
   f^{-1}_1$ is an isomorphism from $f_1(N \rest \delta)$ onto $f_2(N
   \rest \delta)$, mapping $f_1(M \rest \delta)$ onto $f_2(M \rest \delta)$. 
\end{PROOF}

\begin{definition}
\label{inv.10}  Assume $\varphi = \varphi(\bar x,\bar y) \in
\bbL(\tau_T)$ and $N_1 \prec N_2$ are models of $T$ of
cardinality $\lambda$.

\noindent
1) For one-to-one mapping $f$ from $N_2$ to $\lambda$
and $\delta < \lambda$ we define

\begin{equation*}
\begin{array}{clcr}
\text{inv}^\varphi_5(\delta,f,N_1,N_2) = \{\cP \subseteq \cP(\delta):
&\text{there are } 
\bar a_\gamma \in {}^{\ell g(\bar x)}N_2 \text{ satisfying } 
f(\bar a_\gamma) \in {}^{\ell g(\bar x)} \delta \\
  &\text{ for } \gamma < \delta \text{ such that for every} \\
  &\cU \subseteq \delta \text{ the following are equivalent}: \\
  &(i) \quad  {\cU} \in {\cP} \\
  &(ii) \quad \text{for some } \bar b \in {}^{\ell g(\bar y)}N_1 
\text{ we have} \\
  &\quad \quad \gamma < \delta \Rightarrow N_2 \models
  \varphi[\bar a_\gamma,\bar b]^{\text{if}(\gamma \in {\cU})}\}.
\end{array}
\end{equation*}
\mn
2) We let inv$^\varphi_6(N_1,N_2)$ be 
$\langle\text{inv}^\varphi_5(\delta,f,N_1,
N_2):\delta < \lambda\rangle/{\cD}_\lambda$ for some
(equivalently every) $f$ as above.
\end{definition}

\begin{claim}
\label{inv.11}  
1) In Definition \ref{inv.10}(2) we have {\rm inv}$^\varphi_6(N_1,N_2)$
is well defined.

\noindent
2) In Definition \ref{inv.10}, for
$\delta,\lambda,N_1,N_2,\varphi(\bar x,\bar y)$ as there
\mn
\begin{enumerate}
\item[$(a)$]    the set {\rm inv}$^\varphi_5(\delta,f,N_1,N_2)$ 
has cardinality at most $2^{|\delta|}$
\sn
\item[$(b)$]   if $\pi$ is a one-to-one function from $f(N_2)$ into
$\lambda$ mapping $f(N_2) \cap \delta$ onto $\pi(f(N_2)) \cap \delta$
then {\rm inv}$^\varphi_5(\delta,\pi \circ f,N_1,N_2) = 
\text{\rm inv}^\varphi_5(\delta,f,N_1,N_2)$.
\end{enumerate}
\end{claim}

\begin{PROOF}{\ref{inv.11}}  
Easy.
\end{PROOF}

\begin{definition}
\label{inv.12} 
1) For $\varphi = \varphi(\bar x,\bar y)
\in \bbL(\tau_T)$, a model $N$ of $T$ with universe $\lambda,\delta$
a limit ordinal $< \lambda$ and $\kappa < \lambda$ let

\begin{equation*}
\begin{array}{clcr}
\text{inv}^\varphi_{7,\kappa}(\delta,N) = \{{\cP} \subseteq 
{\cP}(\delta):&\text{ we can find } \bar a^i_\gamma \in 
{}^{\ell g(\bar x)}\delta 
\text{ for } \gamma < \delta,i < \kappa \text{ such that} \\
  &\text{ the following conditions on } {\cU} \subseteq \delta 
\text{ are equivalent}: \\
  &(i) \quad {\cU} \in {\cP} \\
  &(ii) \quad \text{for some } \bar b \in {}^{\ell g(\bar y)}N 
\text{ we have}: \\
  &\hskip30pt \text{for every } i < \kappa \text{ large enough for every} \\
  &\hskip30pt  \gamma < \delta \text{ we have }  N \models
  \varphi[\bar a^i_\gamma,\bar b]^{\text{if}(\gamma \in {\cU})}\}.
\end{array}
\end{equation*}
\mn
2) For $\varphi = \varphi(\bar y,\bar x) \in 
\bbL(\tau_T)$ and a model $N$ of $T$ of 
cardinality $\lambda$ let inv$^\varphi_{8,\kappa}(N) = 
\langle\text{inv}^\varphi_7(\delta,N'):\delta < 
\lambda \rangle/{\cD}_\lambda$ for every, equivalently some model
$N'$ isomorphic to $N$ with universe $\lambda$.
\end{definition}

\begin{observation}
\label{inv.13}  
1) inv$^\varphi_{8,\kappa}(N)$ is well defined
for $N \in \text{ EC}_\lambda(T)$ when $|T| + \kappa < \lambda$.

\noindent
2) In Definition \ref{inv.12}(1) we have
$|\text{inv}^\varphi_{7,\kappa}(\delta,N)| \le 2^{|\delta|+\kappa}$.
\end{observation}

\begin{PROOF}{\ref{inv.13}}  
Easy.
\end{PROOF}

\begin{claim}
\label{inv.21}  
Assume $\lambda > |T|$ is regular, $S \subseteq \lambda$ is stationary,
$\varphi = \varphi(\bar x,\bar y)$ and
\mn
\begin{enumerate}
\item[$(a)$]   $\langle N_i:i < \kappa\rangle$ is a $\prec$-increasing
sequence
\sn
\item[$(b)$]   $N_i \in \text{\rm EC}_\lambda(T)$
\sn
\item[$(c)$]   $N = \cup\{N_i:i < \kappa\}$
\sn
\item[$(d)$]   $\bar{\cP} = \langle {\cP}_\alpha:\alpha <
  \lambda\rangle$ where ${\cP}_\alpha \subseteq {\cP}(\alpha)$
\sn
\item[$(e)$]  $f$ is a one-to-one function from $N$ onto $\lambda$
\sn
\item[$(f)$]  for every 
$j < \kappa$ for a club of $\delta$'s below $\lambda$
there are $\bar a^j_\gamma \in N_{j+1} \cap f^{-1}(\delta)$ for $\gamma <
\delta$ satisfying
\begin{enumerate}
\item[$(\alpha)$]   for every $\bar c \in {}^{\ell g(\bar x)}
(N_j)$ there is ${\cU} \in {\cP}_\delta$ such that $\gamma <
\delta \Rightarrow N \models 
\varphi[\bar a^j_\gamma,\bar c]^{\text{\rm if}(\gamma \in {\cU})}$
\sn
\item[$(\beta)$]   for every \, ${\cU} \in {\cP}_\delta$ for some
$\bar b \in {}^{\ell g(\bar y)}(N_\delta)$ we have
$\gamma < \delta\Rightarrow N \models \varphi[\bar a^j_\gamma,
\bar b]^{\text{\rm if}(\gamma \in \cU)}$.
\end{enumerate}
\end{enumerate}
\mn
\Then \, $\{\delta \in S:{\cP}_\delta \in 
\text{\rm inv}^\varphi_{7,\kappa}(\delta,f(N))\} \in {\cD}_\lambda +S$.
\end{claim}

\begin{PROOF}{\ref{inv.21}}
Straight.
\end{PROOF}
\bigskip

\noindent
Now we come to the main two results of this section.
\begin{theorem}
\label{inv.14}  
For some club $E$ of $\lambda^+$, 
if $\delta_1 \ne \delta_2$ belong to $E \cap S^{\lambda^+}_\kappa$
\then \, $M_{\delta_1},M_{\delta_2}$ are not isomorphic, moreover
{\rm inv}$^\varphi_{8,\kappa}(M_{\delta_1}) \ne 
\text{\rm inv}^\varphi_{8,\kappa}(M_{\delta_2})$ \when \,:
\mn
\begin{enumerate}
\item[$\boxtimes$]  $(a) \quad T$ has the strong independence property
(see Definition \ref{0.7}) 
\sn
\item[${{}}$]   $(b) \quad \lambda = \lambda^{<\lambda}$ is 
regular uncountable,
$\lambda > |T|,\lambda > \kappa = \text{\rm cf}(\kappa)$ and
$\lambda^+ = 2^\lambda$
\sn
\item[${{}}$]   $(c) \quad M$ is a saturated model of $T$ of cardinality
$\lambda^+$
\sn
\item[${{}}$]  $(d) \quad \langle M_\alpha:\alpha < \lambda^+\rangle$ is
$\prec$-increasing continuous sequence with union $M$, 

\hskip25pt each of cardinality $\lambda$.
\end{enumerate} 
\end{theorem}

\begin{theorem}
\label{inv.24} Assume $\boxtimes$ of \ref{inv.14}.

\noindent
1)  For some club $E$ of $\lambda^+$, if $\delta_1 < \delta_2 <
\delta_3$ are from $E$ and $\delta_\ell \in S^{\lambda^+}_\lambda$ for
$\ell = 1,2,3$ \then \, $(M_{\delta_2},M_{\delta_1}) \ncong
(M_{\delta_3},M_{\delta_1})$, moreover 
{\rm inv}$^\varphi_6(M_{\delta_2},M_{\delta_1}) \ne 
\text{\rm inv}^\varphi_6(M_{\delta_3},M_{\delta_1})$ for some $\varphi$.

\noindent
2) If $M \prec N_0$ are models of $T$ of cardinality $\lambda$,
\then \, for some elementary extension $N_1 \in 
\text{\rm EC}_\lambda(T)$ of $N_0$ we have $N_1 \prec N_2 \in 
\text{\rm EC}_\lambda(T) \Rightarrow (N_0,M) \not\cong (N_2,M)$.
\end{theorem}

\begin{discussion}  We 
shall below start with $M \in \text{ EC}_\lambda(T)$ and a sequence
$\langle b_i:i < \lambda\rangle$ of distinct members such that
$\langle \varphi(\bar b_i,\bar y):i < \lambda\rangle$ are
independent, and like to find $N,\langle \bar a_i:i <
\lambda\rangle$ such that $M \prec N \in \text{ EC}_\lambda(T)$ and
the $\langle \bar b_i:i < \lambda\rangle$ has a real affect on the relevant
$\varphi$-invariant, in the case of \ref{inv.24}(1) this is
inv$^\varphi_6(M,N)$: for a stationary set of $\delta$'s below 
$\lambda$ it adds
something to the $\delta$-th component in a specific
representation, i.e. assuming $f:N \rightarrow \lambda$ is a
one-to-one function and we deal with $\langle
\text{inv}^\varphi_5(\delta,f,M,N):\delta < \lambda\rangle$; we have
freedom about $\varphi(\bar a_\alpha,\bar b_i)$ and we can assume
$\bar b \in {}^{\ell g(\bar y)}M \backslash 
\{\bar b_i:i < \lambda\} \Rightarrow N \models \neg
\varphi[\bar a_\alpha,\bar b]$.

But the relevant ${\cP}_\delta$ is influenced not just  by say
$\langle \bar b_i:i \in [\delta,2^{|\delta|})\rangle$ but also by
later $\bar b_i$'s (and earlier $\bar b_i$).  
To control this we use below $\langle \bar
a_\alpha:\alpha < \lambda\rangle,S,E$ such that we deal with different
$\delta \in S$ in an independent way to large extent; 
this is the reason for choosing the $C^*_\alpha$'s.
\end{discussion}

\begin{PROOF}{\ref{inv.14}}
\underline{Proof of \ref{inv.14}}  By the proof of \cite[\S2]{Sh:877} 
\wilog \, $\lambda$ is strongly
inaccessible.  Choose $\theta \in \text{ Reg } \cap \lambda \backslash
\{\aleph_0\}$, will be needed when we generalize the proof in \S2.

Let $\langle {\cU}_i:i < \kappa\rangle$ be a $\subseteq$-increasing
sequence of subsets of $\lambda$ such that the set
${\cU}^-_i = {\cU}_i \backslash \cup\{{\cU}_j:j < i\}$ has cardinality
$\lambda$ for each $i < \kappa$ and
let $\cU_\kappa = \cup\{\cU_i:i < \kappa\}$.  
Let $\varphi(\bar x,\bar y) \in \bbL(\tau_T)$
have the strong independence property, see Definition \ref{0.7}.  We
can choose $\langle C^*_\alpha:\alpha < \lambda^+\rangle$ such that
$C^*_\alpha \subseteq \text{ nacc}(\alpha)$, otp$(C^*_\alpha) \le
\kappa,\beta \in C^*_\alpha \Rightarrow C^*_\beta = C^*_\alpha \cap
\beta$ and $\lambda|\alpha \wedge \text{ cf}(\alpha) = \kappa
\Rightarrow \alpha = \sup(C^*_\alpha)$ and cf$(\alpha) \ne \kappa
\Rightarrow \text{ otp}(C^*_\alpha) < \kappa$.

[See \cite{Sh:237e} but for 
completeness we show this; by induction on $\alpha < \lambda^+$
we choose $\langle
C^*_\varepsilon:\varepsilon < \lambda \alpha\rangle$ such that:
\mn
\begin{enumerate}
\item[$(a)$]  the relevant demand holds
\sn
\item[$(b)$]  if $\alpha = \beta +1,C \subseteq \lambda \beta,(\forall
i \in C)(C^*_i = C \cap i)$ and otp$(C) < \kappa$ then for some $i \in
(\lambda \beta,\lambda \alpha)$ we have $C^*_{i+1} = C$.
\end{enumerate}
\mn
As $\lambda = \lambda^{< \kappa}$ this is easy but we elaborate.  For
$\alpha = 0$ trivial for $\alpha$ limit obvious.  Assume $\alpha =
\beta +1$ let $\alpha_* = \lambda \alpha,\beta_* = \lambda \beta$ and
$\langle C^*_i:i < \beta_*\rangle$ has been defined.

First, we choose $C^*_{\beta_i}$.  If cf$(\beta_*) \ne \kappa$ let
$C^*_{\beta_*} = \emptyset$, so assume cf$(\beta_*) = \kappa$ then
necessarily cf$(\beta) = \kappa$.

Let $\langle \alpha_\varepsilon:\varepsilon <
\kappa \rangle $ be increasing with limit $\beta$, and choose
$\beta_\varepsilon \in [\lambda \alpha_\varepsilon,\lambda
\alpha_\varepsilon + \lambda)$ by induction on $\varepsilon < \kappa$
such that $C^*_{\beta_\varepsilon} = \{\beta_\zeta:\zeta <
\varepsilon\}$.  

Lastly, let $C^*_{\lambda \alpha} :=
\{\beta_\varepsilon:\varepsilon < \kappa\}$.  So $C^*_{\beta_*}$ has
been defined in any case.

Now let $\cC_\alpha =
\{C \subseteq \lambda \alpha:\text{otp}(C) < \kappa$ and $\beta \in C
\Rightarrow C^*_\beta = C \cap \beta\}$, so $|\cC| \le \lambda$, also
$\emptyset \in \cC$, so let $\langle C^*_{\lambda \alpha +i}:i \in
(0,\lambda)\rangle$ list $\cC_\alpha$ possibly with repetitions.
So we have defined
$\langle C^*_\varepsilon:\varepsilon < \lambda \alpha\rangle$, so 
have carried the induction.]

Let $S_* = \{\mu:\mu = \beth_{\alpha + \omega}$ for some $\alpha <
\lambda\}$.  
Let $E_*,\langle C_\alpha:\alpha < \lambda\rangle$ be such that:
\mn
\begin{enumerate}
\item[$\circledast_1$]   $(a) \quad C_\alpha \subseteq \alpha \cap S_*$
\sn
\item[${{}}$]  $(b) \quad \beta \in C_\alpha \Rightarrow C_\beta =
C_\alpha \cap \beta$
\sn
\item[${{}}$]   $(c) \quad$ otp$(C_\alpha) \le \theta$
\sn
\item[${{}}$]   $(d) \quad E_*$ is the  club $\{\delta <
\lambda:\theta < \delta = \beth_\delta\}$ of $\lambda$
\sn
\item[${{}}$]   $(e) \quad$ 
otp$(C_\alpha) = \theta$ iff $\alpha \in E_* \cap S^\lambda_\theta$
\sn
\item[${{}}$]  $(f) \quad$ if $\alpha \in S := E_* \cap S^\lambda_\theta$
then $\alpha = \sup(C_\alpha)$
\sn
\item[${{}}$]  $(g) \quad$ if $\alpha \in E_*$ and $i < \kappa$ then
$|\alpha \cap \cU^-_i| = |\alpha|$.
\end{enumerate}
\mn
[Why can we choose?  By induction on the cardinal $\chi \ge \aleph_0$
we choose $\langle C_\alpha:\alpha < \beth_\chi\rangle$ and $E_\chi =
E_* \cap \beth_\chi$ such that the
relevant demands hold and: if $\chi = 2^{\chi_1}$ and $C$ is a subset
of $S_* \cap \beth_{\chi_1}$ of order type $< \theta$ satisfying
$\alpha \in C \Rightarrow C_\alpha = C \cap \alpha$ \then \, for some
$\alpha \in S_* \cap (\beth_{\chi_1},\beth_\chi)$ we have $C_\alpha =
C$.  Why this extra induction hypothesis help?  As arriving to $\alpha
\in S$ so $\alpha = \beth_\alpha$ let $\langle \chi_i:i <\theta
\rangle$ be an increasing sequence of cardinals with limit
$\beth_\alpha = \alpha$ and we choose $\alpha_i \in
(\beth_{\chi_i},\beth_{2^{\chi_i}}) \cap S_*$ by induction on $i <
\theta$ such that $C_{\alpha,i} = \{\alpha_j:j < i\}$ and the let
$C_\chi = \{\alpha_i:i < \theta\}$.]

We shall prove that
\mn 
\begin{enumerate}
\item[$\circledast_2$]   if $\boxdot_2$ below holds, then 
 there is a $\beta$ such that $\odot_2$ holds \underline{where}:
\begin{enumerate}
\item[$\boxdot_2$]  $(a) \quad \alpha < \lambda^+,i < \kappa$
\sn
\item[${{}}$]   $(b) \quad f$ is a one-to-one function from
$M_\alpha$ into $\cU'_i = \cup\{\cU_j:j < i\}$
\sn
\item[${{}}$]   $(c) \quad E \subseteq E_*$ is a club of 
$\lambda$ such that $\delta \in E \Rightarrow
f(M_\alpha) \restriction \delta \prec f(M_\alpha)$
\sn
\item[${{}}$]    $(d) \quad \bar{\cP} = \langle {\cP}_\delta:
\delta \in S\rangle$
\sn
\item[${{}}$]    $(e) \quad {\cP}_\delta \subseteq 
{\cP}(\delta)$ and $\emptyset \in {\cP}_\delta$ and ${\cP}_\delta
\subseteq \bigcup\limits_{\ell \le 2} {\cP}^{*,\ell}_\delta$
where\footnote{note that $\cP^{*,1}_{\delta,i},\cP^{*,2}_{\delta,i}$
are the families of sets 
we like to ignore as they are influenced by our choices
for $\delta_1 \in S \backslash \{\delta\}$, so we work to have them
families of bounded subsets of $\delta$.} 
\sn
\item[${{}}$]  \hskip20pt $(\alpha) \quad {\cP}^{*,0}_\delta = 
\{A \subseteq \delta:\sup(A) = \delta$ and $A
\subseteq \cup\{[\mu,2^\mu):\mu \in C_\delta\}\}$,
\sn
\item[${{}}$]  \hskip20pt $(\beta) \quad {\cP}^{*,1}_\delta =
\cup\{{\cP}^{*,0}_{\delta_1}:\delta_1 \in S \cap \delta\}$,
\sn
\item[${{}}$]  \hskip20pt $(\gamma) \quad {\cP}^{*,2}_\delta = \{A
\subseteq \delta$\,: for some $\delta_1 \in \lambda \backslash (\delta +1)$
we have

\hskip40pt $A \subseteq 
\cup\{[\partial,2^\partial):\partial \in C_{\delta_1} \cap \delta\}\}$
\sn
\item[${{}}$]   $(f) \quad$ if $\delta_1 < \delta_2$ are from $S$ 
\then 
\sn
\item[${{}}$]  \hskip20pt $(\alpha) \quad A \in {\cP}_{\delta_1} 
\Rightarrow A \in {\cP}_{\delta_2}$ 
\sn
\item[${{}}$]  \hskip20pt $(\beta) \quad A \in {\cP}_{\delta_2} 
\Rightarrow A \cap \delta_1 \in {\cP}_{\delta_1}$,
\sn
\item[${{}}$]  \hskip20pt $(\gamma) \quad$ for any $\delta \in S$ the family
${\cP}^{*,1}_\delta \cup {\cP}^{*,2}_\delta$ is a set of bounded 

\hskip45pt subsets of $\delta$; (this follows)
\sn
\item[${{}}$]   $(g) \quad \bar b_{\delta,{\cU}} \in M_\alpha$ for $\delta
\in E \cap S,\cU \in \cP_\delta$ are such that 

\hskip40pt $\bar b_{\delta_1,{\cU}_1}
= \bar b_{\delta_2,{\cU}_2} \wedge \cU_1 \in \cP_{\delta_1} \wedge \cU_2
\in \cP_{\delta_2} \Rightarrow \delta_1 = \delta_2 \wedge {\cU}_1 = {\cU}_2$
\sn
\item[$\odot_2$]  $(\alpha) \quad \beta \in (\alpha,\lambda^+)$
\sn
\item[${{}}$]    $(\beta) \quad$ there are $\bar a_\gamma \subseteq M_\beta$  
for $\gamma < \delta$ such that for a club of $\delta \in E$, if

\hskip25pt $\delta \in S$ 
then the following conditions on ${\cU} \subseteq \delta$ are
equivalent:
\sn
\item[${{}}$]  \hskip45pt $(i) \quad {\cU} \in {\cP}_\delta$
\sn
\item[${{}}$]  \hskip45pt  $(ii) \quad$ for some 
$\bar b \in {}^{\ell g(\bar y)}M_\alpha$ we have: for every $\gamma < \delta$,

\hskip75pt $M_\beta \models \varphi[\bar a_\gamma,\bar b]$
iff $\gamma \in {\cU}$
\sn
\item[${{}}$]  \hskip45pt $(iii) \quad$  clause (ii) 
holds for $\bar b = \bar b_{\delta,{\cU}}$ and $\cU \in \cP_{\delta,i}$
\end{enumerate}
\end{enumerate}
\mn
[Why?  For each $\delta \in E \cap S$ let
$\langle {\cU}_{\delta,\varepsilon}:\varepsilon < |{\cP}_{\delta,i}| \le
2^{|\delta|} \rangle$ list ${\cP}_\delta$ and let 
$\bar b_{\delta,\varepsilon} := \bar b_{\delta,{\cU}_{\delta,\varepsilon}}$.

Let

\begin{equation*}
\begin{array}{clcr}
\Gamma = \{\varphi(\bar x_\gamma,\bar b_{\delta,\varepsilon})^{\text{if}(\gamma
  \in {\cU}_{\delta,\varepsilon})}:&\gamma < \lambda,
\delta \in E \text{ and } \varepsilon <   |{\cP}_{\delta,i}|\} \\
  &\cup\{\neg \varphi(\bar x_\gamma,\bar b):\gamma < \lambda,
\bar b \in {}^{\ell g(\bar y)}(M_\alpha) \text{ and for no} \\
  & \delta \in E,\varepsilon < |{\cP}_\delta| \text{ do we
have } \bar b = \bar b_{\delta,\varepsilon}\}.
\end{array}
\end{equation*}

\mn
As $\varphi(\bar x,\bar y)$ has the strong independence property,
recalling that by clause (g) of $\boxdot_2$ the sequence $\langle \bar
b_{\delta,\varepsilon}:\delta \in E \cap S$ and 
$\varepsilon < |\cP_\delta|\rangle$ is
with no repetitions, clearly
$\Gamma$ is finitely satisfiable in $M_\alpha$, but $M$ is 
$\lambda^+$-saturated,
$M_\alpha \prec M$ and $|\Gamma| = \lambda$ hence we can find
$\bar a_\gamma \in {}^{\ell g(\bar x)}M$ 
for $\gamma < \lambda$ such that the assignment
$\bar x_\gamma \mapsto \bar a_\gamma \, (\gamma < \lambda)$ satisfies $\Gamma$
in $M$.  Lastly, choose $\beta \in (\alpha,\lambda^+)$ such that
$\{\bar a_\gamma:\gamma < \lambda\} \subseteq M_\beta$.

Now check recalling $\emptyset \in \cP_\delta$ for $\delta \in S$.]

Note
\mn
\begin{enumerate}
\item[$\odot_3$]   in $\odot_2$ if
$h$ is a one-to-one mapping from $M_\beta$ into $\cU$ extending $f$
\then \, for some club $E$ of $\lambda$ if for every $\delta \in S \cap E$
we have $(\forall \gamma < \lambda)(\gamma < \delta \rightarrow h(\bar
a_\gamma) \in {}^{\ell g(\bar x)}\delta)$ and so for every $\cU
\subseteq \delta$ the conditions $(i),(ii),(iii)$ from $\odot_2$ are
equivalent.
\end{enumerate}
\mn
Next we can choose $\bar f$ such that
\mn
\begin{enumerate}
\item[$\circledast_3$]  $(a) \quad \bar f = \langle f_\alpha:\alpha
  < \lambda^+\rangle$
\sn
\item[${{}}$]  $(b) \quad f_\alpha$ is a one-to-one function from
  $M_\alpha$ into $\cU_{\text{otp}(C^*_\alpha)}$
\sn
\item[${{}}$]  $(c) \quad$ if $\alpha \in
C^*_\beta$ then $f_\alpha \subseteq f_\beta$.
\end{enumerate}
\mn
Now
\mn
\begin{enumerate}
\item[$\circledast_4$]  for every $\alpha < \lambda^+$ there is
$\bar{\cP}^\alpha = \langle {\cP}^\alpha_\varepsilon:
\varepsilon \in S\rangle$ such that
\begin{enumerate}
\item[$(i)$]  ${\cP}^\alpha_\varepsilon \subseteq 
{\cP}(\varepsilon)$ are as in $\boxdot_2(e),(f)$ above
\sn
\item[$(ii)$]  for every $\beta \le \alpha$, for a club of
$\delta$'s from $S$ we have ${\cP}^\alpha_\delta \notin 
\text{ inv}^\varphi_{7,\kappa}(\delta,f_\beta(M_\beta))$. 
\end{enumerate}
\end{enumerate}
\mn
[Why?  For every $\beta \le \alpha$ and $\delta \in (\kappa,\lambda)$
we have inv$^\varphi_{7,\kappa}(\delta,f_\beta(M_\beta))$ is a subset of
${\cP}({\cP}(\delta))$ of cardinality $\le 2^{|\delta|}$.  As
the number of $\beta$'s is $\le \lambda$, by diagonalization we can do this:
 let $\alpha +1 = \bigcup\limits_{\varepsilon < \lambda} 
u_\varepsilon$ and $u_\varepsilon \in [\alpha+1]^{< \lambda}$ 
increasing continuous for $\varepsilon < \lambda$;
moreover, $\alpha < \lambda \Rightarrow u_\varepsilon = \alpha$ and
$\alpha \ge \lambda \Rightarrow 
u_\varepsilon \cap \lambda \subseteq \varepsilon$ and 
$|u_\varepsilon| \le |\varepsilon|$.   
By induction on $\varepsilon \in (\kappa,\lambda) \cap S$ choose 
${\cP}^\alpha_\varepsilon \subseteq \bigcup\limits_{\ell < 3}
{\cP}^{*,\ell}_{\alpha_\varepsilon}$ 
which includes $\cup\{{\cP}^\alpha_\zeta:\zeta \in
u_\varepsilon \cap S\} \cup {\cP}^{*,2}_\varepsilon$ and satisfies
${\cP}^{*,0}_\varepsilon \cap {\cP}^\alpha_\varepsilon \in 
{\cP}({\cP}^{*,0}_{\varepsilon,i})) \backslash \{\cP \cap
\cP^{*,0}_\varepsilon:\cP \in \text{ inv}^\varphi_{7,\kappa}
(\varepsilon,f_\beta(M_\beta)),\beta \in u_\varepsilon\}$.  Note that
for each $\beta \le \alpha$ the set $\{\varepsilon < \lambda:\beta \in
u_\varepsilon\}$ contains an end-segment of $\lambda$ hence a club of
$\lambda$ as required.]

Now choose pairwise distinct $\bar b_{\delta,{\cU}} \in {}^{\ell
g(\bar y)}(M_0)$ for $\delta \in E_*,{\cU} \in {\cP}^{*,0}_\delta$
\mn
\begin{enumerate}
\item[$\circledast_5$]  for every $\alpha_* \le \alpha < \lambda^+$
for some $\beta \in (\alpha,\lambda^+)$ and $\bar a_\gamma \in 
{}^{\ell g(\bar x)} M_\beta$ for $\gamma < \lambda$ the condition 
in clause $(\gamma)$ of $\odot_2$ holds with $\bar{\cP}^{\alpha_*}$
here standing for $\bar{\cP}$ there and the $\bar b_{\delta,{\cU}}$
chosen above.
\end{enumerate}
\mn
[Why?  By $\circledast_2$.]
\mn
\begin{enumerate}
\item[$\circledast_6$]   let $E = \{\delta < \lambda^+:\delta$ is a
  limit ordinal such that for every $\alpha_* \le \alpha < \delta$ there
  is $\beta < \delta$ as in $\circledast_5\}$.
\end{enumerate}
\mn
Clearly $E$ is a club of $\lambda^+$.
\mn
\begin{enumerate}
\item[$\circledast_7$]   if $\delta_1 < \delta_2$ are from $E \cap 
S^{\lambda^+}_\kappa$ then $M_{\delta_1},M_{\delta_2}$ are not
isomorphic.
\end{enumerate}
\mn
[Why?   Let $\alpha_* = \text{ min}(C^*_{\delta_2} \backslash \delta_1)$.
We consider $\bar{\cP}^{\alpha_*}$ which is from $\circledast_4$.
On the one hand $\{\varepsilon < \lambda:{\cP}^{\alpha_*}_\varepsilon \notin
\text{ inv}^\varphi_{7,\kappa}(\varepsilon,f_{\delta_1}
(M_{\delta_1}))\}$  contains a club by $\circledast_4(ii)$.  Note that
$\langle f_\alpha:\alpha \in C^*_{\delta_2} \backslash
\delta_1\rangle$ is $\subseteq$-increasing sequence of functions 
with union $f_{\delta_2}$.

\noindent
On the other hand choose an increasing $\langle \alpha_i:i < \kappa\rangle$
 with limit $\delta_2$ satisfying $\alpha_0 = 0,\alpha_1 = \delta_1$ such that
$(\alpha_*,\alpha_{1+i},\alpha_{1+i+1})$ are like 
$(\alpha_*,\alpha,\beta)$ in $\circledast_5$ for each $i < \kappa$ and
 $i \in (1,\kappa) \Rightarrow \alpha_i \in C^*_{\delta_2}$.  Now
 by \ref{inv.21}, $\{\varepsilon < \lambda:{\cP}^{\alpha_*}_\varepsilon \in 
\text{ inv}^\varphi_{7,\kappa}
(\varepsilon,f_{\delta_2}(M_{\delta_2}))\}$  contains a club.
Hence by the last sentence and the end of the previous paragraph
$M_{\delta_1} \ncong M_{\delta_2}$ as required.]

So we are done.  
\end{PROOF}

\begin{remark}  We can avoid using $C^*_\delta$ and also $C_\delta$
(e.g. using $A \in \cP^{*,0}_\delta \Rightarrow \text{ otp}(A)=\delta$)
but seems less transparent.
\end{remark}

\begin{PROOF}{\ref{inv.24}}
\underline{Proof of \ref{inv.24}}  Similar but easier (for $\lambda$
regular not strong limit (but $2^\lambda > 2^{< \lambda}$) also easy),
or see the proof of \ref{3e.14}.  
\end{PROOF}

\section{Independent $T$} 

We would like to do something similar to \S1, but our control on the
relevant family of subsets of $\mu$ is less tight.  We control it to
some extent by using the completion of a free Boolean algebra.

\begin{context}
\label{2d.1}   $T$ a complete first order theory,
$\varphi(x,\bar y)$ has the independence property (of course the
existence of such $\varphi$ follows from the strong independence
property but is weaker).

We continue \cite[2.1-2.12]{Sh:877}, but we do not rely on it.
\end{context}

\begin{definition}
\label{2d.4}  For a set $I$ let
\mn
\begin{enumerate}
\item[$(a)$]   $\bbB = \bbB_I$ be the Boolean Algebra generated
  by $\langle e_t:t \in I\rangle$ freely,
\sn
\item[$(b)$]  $\bbB^c_I$ is the completion of $\bbB$
\sn
\item[$(c)$]   for $J \subseteq I$ let $\bbB^c_{I,J}$ be the
complete subalgebra of $\bbB^c_I$ generated by $\{e_s:s \in J\}$
\sn
\item[$(d)$]   let uf$(\bbB^c_I)$ be the set of ultrafilters on $I$.
\end{enumerate}
\end{definition}

\begin{claim}
\label{2d.7}  Assume
\mn
\begin{enumerate}
\item[$\circledast$]   $(a) \quad M \models T$
\sn
\item[${{}}$]   $(b) \quad \bar b_t \in {}^{\ell g(\bar y)}M$ for
$t \in I$
\sn
\item[${{}}$]   $(c) \quad \langle \varphi(x,\bar b_t):t \in
I\rangle$ is an independent sequence of formulas.
\end{enumerate}
\mn
\Then \, there is a function $F$ from ${}^{\ell g(\bar y)}M$ to $\bbB
= \bbB^c_I$ such that 
\mn
\begin{enumerate}
\item[$(\alpha)$]   $F(\bar b_t) = e_t$
\sn
\item[$(\beta)$]   for every ultrafilter $D$ of $\bbB$ there is
$p = p_D = p_{F,D} \in \bold S_\varphi(M)$, in fact, a unique one, 
such that for every $\bar b \in {}^{\ell g(\bar y)}M$ we have 
$\varphi(x,\bar b) \in p \Leftrightarrow F(\bar b) \in D$.
\end{enumerate}
\end{claim}

\begin{remark}  1) Note that the mapping $D \mapsto p_D$ is not
necessarily one to one, but $D_1 \cap \{e_t:t \in I\} \ne D_2 \cap
\{e_t:t \in I\}\Rightarrow p_{D_1} \ne p_{D_2}$.

\noindent
2) If $I = I_1 \cup I_2,I_1 \cap I_2 = \emptyset$ and $|I_2| =
   |I_1|^{\aleph_0}$ \then \, we can find a mapping $F$ from ${}^{\ell
g(\bar y)}M$ onto (not just into) $\bbB = \bbB^c_{I_1}$ such that
clause $(\alpha),(\beta)$ are satisfied.
\end{remark}

\begin{PROOF}{\ref{2d.7}}
Clearly ${\cP}(M)$ is a Boolean algebra and $\{\varphi(M,\bar
b_t):t \in M\}$ generates freely a subalgebra of ${\cP}(M)$ 
which we call $\bbB'$.  So there is a homomorphism 
$h$ from $\bbB'$ into $\bbB$ mapping $\varphi(M,\bar b_t)$ to 
$e_t$ (moreover $h$ is unique and is an isomorphism from $\bbB'$ onto 
$\bbB_I \subseteq \bbB^c_I$).  So $h$ is a
homomorphism from $\bbB' \subseteq {\cP}(M)$ into $\bbB^c_I$,
which is a complete Boolean algebra hence there is a homomorphism
$h^+$ from the Boolean algebra ${\cP}(M)$ into $\bbB^c$ extending $h$.  

Lastly, define $F:{}^{\ell g(\bar y)}M \rightarrow \bbB^c$ by
$F(\bar b) = h^+(\varphi(M,\bar b))$.  Now check.
\end{PROOF}

\begin{conclusion}
\label{2d.14}  Assume $\circledast$ from \ref{2d.7} and
\mn
\begin{enumerate} 
\item[$\boxdot$]  $(a) \quad I = \lambda$ is regular uncountable
\sn
\item[${{}}$]   $(b) \quad |M| \subseteq {\cU} \subseteq \lambda$
\sn
\item[${{}}$]   $(c) \quad D_\alpha$ is an ultrafilter of $\bbB^c_I$
for $\alpha < \lambda$
\sn
\item[${{}}$]   $(d) \quad {\cU} \backslash |M|$ is unbounded in $\lambda$.
\end{enumerate}
\mn
\Then \, we can find $\langle a_\alpha:\alpha < \lambda\rangle$ and
$N$ such that
\mn
\begin{enumerate}
\item[$(\alpha)$]  $M \prec N$
\sn
\item[$(\beta)$]  $|N| \subseteq {\cU}$
\sn
\item[$(\gamma)$]  $a_\alpha \in N$ for $\alpha < \lambda$
\sn
\item[$(\delta)$]   $a_\alpha$ realizes $p_{D_\alpha} \in 
\bold S_\varphi(M)$.
\end{enumerate}
\end{conclusion}

\begin{remark}
Conclusion \ref{2d.14} is easy but intended to clarify how we shall use
the ultrafilters, so is quoted toward the end of the section.
\end{remark}

\begin{PROOF}{\ref{2d.14}}
 Should be clear.
\end{PROOF}

\begin{discussion}
\label{2d.21}  Note that compared to \S1 instead $\bar
x,\bar y,\bar a_\alpha,\bar b_\beta$ we have $x,\bar y,a_\alpha,\bar b_\beta$. 
 Compared to \S1, we have less control over
$\{\text{tp}(a,M,N):a \in N\}$.  There, for the sequences $\bar b$ of $M$ which
are not among $\{\bar b_\gamma:\gamma < \lambda\}$, 
we can demand $N \models \neg \varphi[\bar a_\gamma,\bar b]$ for
$\gamma < \lambda$ so tp$_\varphi(\bar a_\gamma,M,N)$ can be clearly read.
Here the complete Boolean Algebra $\Bbb B^c_I$ is helping, a small
price is that we need $\theta > \aleph_0$.

In order to try to keep track of what is going on we shall use only
 tp$(a_\gamma,M,N)$ of the form $p_D$ for ultrafilter $D$ on 
$\bbB^c_I$.  Further, we better have, e.g. a nice function $\pi$
from ${}^\lambda 2$ to uf$(\bbB^c_I)$ such that $(e_\alpha \in
\pi(\eta)) \Leftrightarrow \eta(\alpha)=1$.

A possible approach is: we define $\langle M_{\eta,u}:\eta \in {\cT}
\subseteq \text{ des}(\lambda),u \in {\cP}(n_\eta)\rangle$ as in
\cite[\S3]{Sh:668} and we define $D_\eta \in \text{ uf}(\bbB^c \cap M)$
such that $\alpha \in M_\eta \cap \lambda \Rightarrow
[e_\eta^{\eta(\alpha)} \in D_\eta]$ and $\bigcup\limits_{\eta} \bbD_\eta \in 
\text{ uf}(\bbB^c)$.  

We need some continuity so each ``$e \in D_\eta$" \, $(e \in \bbB^c)$
depends on $\eta \restriction u_e$ for some ``small" $u_e \subseteq
\lambda$.
\end{discussion}

\begin{theorem}
\label{3e.7} 
In Theorem \ref{inv.14} it suffices to
assume $\boxtimes'$ which means clauses (b),(c),(d) of $\boxtimes$ and
\mn
\begin{enumerate}
\item[$(a)'$]  $T$ has the independence property.
\end{enumerate}
\end{theorem}

\begin{theorem}
\label{3e.14}   In Theorem \ref{inv.24} it suffices
to assume $\boxtimes'$ of \ref{3e.7}.
\end{theorem}

\begin{PROOF}{\ref{3e.7}}
\underline{Proof of \ref{3e.7}}  
Just combine the proofs of \ref{inv.14} from \S1 and \ref{3e.14} below.
\end{PROOF}

\begin{PROOF}{\ref{3e.14}}
\underline{Proof of \ref{3e.14}}  
As in the proof of \ref{inv.14} we
can assume $\lambda$ is strongly inaccessible though the proof is just
easier otherwise.  We let
\mn
\begin{enumerate}
\item[$\circledast_1$]   $(a) \quad E_* = \{\delta < \lambda:\delta =
\beth_\delta\}$, a club of $\lambda$
\sn
\item[${{}}$]  $(b) \quad S_* = \{\beth_{\alpha + \omega}:\alpha <
\lambda\}$
\sn
\item[${{}}$]  $(c) \quad$ choose a regular uncountable $\theta <\lambda$
\end{enumerate}
\mn
and let 
\mn
\begin{enumerate}
\item[$\circledast_2$]   $(a) \quad S = \{\delta \in E_*\,$:\,cf$(\delta) =
\theta\} = S^\lambda_\theta \cap E_*$ and
\sn
\item[${{}}$]  $(b) \quad$ let $\bar C$ be as in $\circledast_1$ of
the proof of \ref{inv.14}, in particular $\bar C =$

\hskip25pt $\langle C_\alpha:\alpha < \lambda\rangle,
C_\alpha \subseteq S_*,\text{ otp}(C_\alpha) \le \theta,\alpha \in
C_\beta \Rightarrow C_\alpha = C_\beta \cap \alpha$

\hskip25pt  and $\alpha \in S \Leftrightarrow \alpha = \sup(C_\alpha)
\Leftrightarrow
\text{ otp}(C_\alpha) = \theta$ 

\hskip25pt and $\alpha \in S \Rightarrow \sup(C_\alpha) = \alpha$
\sn
\item[${{}}$]  $(c) \quad$ for $\mu \in S$ let $A^*_\mu =
\cup\{[\chi,2^\chi]:\chi \in C_\mu\}$.
\end{enumerate}
\mn
Let $\bbD_*$ be an ultrafilter of $\bbB^c_\lambda$ such that $e_\alpha
\notin \bbD_*$ for $\alpha < \lambda$.

Now for $\eta \in {}^\lambda 2$ we choose $\bbD_\eta$ such that
\mn
\begin{enumerate}
\item[$\circledast_3$]  $(a) \quad \bbD_\eta$ is an ultrafilter
of $\bbB^c_\lambda$
\sn
\item[${{}}$]  $(b) \quad$ if $e \in \bbD_* \subseteq \bbB^c_\lambda$ 
belongs to $\bbB^c_{\lambda,\eta^{-1}\{0\}}$ (see \ref{2d.4}, the 
completion of  

\hskip25pt the subalgebra of $\bbB^c_\lambda$ 
generated by $\{e_\alpha:\eta(\alpha)=0\}$) then $e \in \bbD_\eta$.
\sn
\item[${{}}$]   $(c) \quad$ if $\alpha < \lambda$ and $\eta(\alpha)=1$
then $e_\alpha \in \Bbb D_\eta$.
\end{enumerate}
\mn
So 
\mn
\begin{enumerate}
\item[$\circledast_4$]  $(a) \quad$ if $\eta \in {}^\lambda 2$ is constantly
zero then $\bbD_\eta  = \bbD_*$
\sn
\item[${{}}$]   $(b) \quad e_\alpha \in \bbD_\eta \Leftrightarrow
\eta(\alpha) = 1$ for $\alpha < \lambda,\eta \in {}^\lambda 2$.
\end{enumerate}
\mn
Now let $\bar \eta = \langle \eta_\varepsilon:\varepsilon <
\lambda\rangle$ be a sequence of members of ${}^\lambda 2$ 
and below we shall be interested mainly
in the case $\alpha = \mu \in S$.  

Define  
\mn
\begin{enumerate}
\item[$\circledast_5$]   for $e \in \bbB^c_\lambda$ and $\alpha
\le \lambda$ we let $Y^\alpha_{\bar\eta,e} := 
\{\varepsilon < \alpha:e \in \bbD_{\eta_\varepsilon}\}$
\sn
\item[$\circledast_6$]   ${\cP}_{\bar\eta,\alpha} := 
\{Y^\alpha_{\bar\eta,e}:e \in \bbB^c_\lambda \bigr\}$.
\end{enumerate}
\mn
Now what can we say on ${\cP}_{\bar \eta,\mu}$ for
$\mu \in S$ ? 
\mn
As we can consider $e \in \{e_\alpha:\alpha \in [\mu,2^\mu)\}$, clearly
\mn
\begin{enumerate}
\item[$\circledast_7$]   $\bigl\{ \{\varepsilon <
\mu:\eta_\varepsilon(\alpha) = 1\}:\alpha \in [\mu,2^\mu) \bigl\} \subseteq
{\cP}_{\bar\eta,\mu} \subseteq {\cP}(\mu)$.
\end{enumerate}
\mn
This may be looked at as a lower bound of ${\cP}_{\bar\eta,\mu}$.
Naturally we try to get also an ``upper bound" to 
${\cP}_{\bar\eta,\mu}$; now note
\mn
\begin{enumerate}
\item[$\circledast_8$]  if $e \in \bbB^c_\lambda$ then
$Y^\mu_{\bar\eta,-e} = \mu \backslash Y^\mu_{\bar\eta,e}$.
\end{enumerate}
\mn
Now define (recalling $A^*_\mu$ is from $\circledast_2(c))$
\mn
\begin{enumerate}
\item[$\circledast_9$]  $\Xi$ is the set of $\bar\eta$ of
the form $\langle \eta_\varepsilon:\varepsilon < \lambda\rangle$ such that:
\begin{enumerate}
\item[$(a)$]  $\eta_\varepsilon \in {}^\lambda 2$ for every
$\varepsilon < \lambda$,
\sn
\item[$(b)$]  if $\eta_\varepsilon(\alpha) =1$ \then
\, $(\exists \mu \in S)[\mu \le \alpha < 2^\mu \wedge \varepsilon \in
A^*_\mu]$, 
\sn
\item[$(c)$]  if $\mu \in S$ and $u \subseteq [\mu,2^\mu)$ is
countable \then \, $\{\varepsilon \in A^*_\mu$: if $\alpha \in u$ then
$\eta_\varepsilon(\alpha) = 0\}$ is of cardinality $\mu$.
\end{enumerate}
\end{enumerate}
\mn
Also (by our knowledge of the completion of a free Boolean algebra, 
$\bbB^c_\lambda$ satisfies the c.c.c.)
for every $e \in \bbB^c_\lambda$ we can choose $u_e$ such that:
\mn
\begin{enumerate}
\item[$\boxplus_1$]   $(a) \quad u_e \subseteq \lambda$ is countable
\sn
\item[${{}}$]   $(b) \quad e \in \bbB^c_{\lambda,u_e}$.
\end{enumerate}
\mn
So by clause (b) of $\circledast_3$ clearly
\mn
\begin{enumerate}
\item[$\boxplus_2$]  if $\bar \eta \in \Xi,e \in \bbB^c_\lambda,
\varepsilon < \mu \in S$ and $u_e \subseteq
\eta^{-1}_\varepsilon\{0\}$ \then \, $e \in \bbD_{\eta_\varepsilon}
\Leftrightarrow e \in \bbD_*$
\end{enumerate}
\mn
hence
\mn
\begin{enumerate}
\item[$\boxplus_3$]   if $\bar\eta \in \Xi,e \in \bbB^c_\lambda 
\cap \bbD_*$ and $\mu \in S$ then $Y^\mu_{\bar\eta,e} \supseteq
\{\varepsilon < \mu:u_e \subseteq \eta^{-1}_\varepsilon\{0\}\}$.
\end{enumerate}
\mn
Next
\mn
\begin{enumerate}
\item[$\boxplus_4$]   for $\bar\eta \in \Xi$
\begin{enumerate}
\item[$(a)$]   let ${\cD}_{\bar\eta,\mu}$ be the filter on
$\mu$ generated by $\{A^*_\mu\} \cup \{\{\varepsilon < \mu:u \subseteq
\eta^{-1}_\varepsilon\{0\}$ 
and $\varepsilon > \zeta\}:\zeta < \mu$
and $u \subseteq \lambda$ is countable$\}$
\sn
\item[$(b)$]   let ${\cI}_{\bar\eta,\mu}$ be the dual ideal. 
\end{enumerate}
\end{enumerate}
\mn
Clearly
\mn
\begin{enumerate}
\item[$\boxplus_5$]   $(a) \quad$ if $\bar\eta \in \Xi,\mu \in S$ and
$\alpha \in \lambda \backslash [\mu,2^\mu)$ then $\{\varepsilon <
\mu:\eta_\varepsilon(\alpha) \ne 0\}$ is

\hskip25pt  a bounded subset of $\mu$,
\sn
\item[${{}}$]   $(b) \quad$ if $\bar\eta \in \Xi$ and $\mu \in
S$ then ${\cD}_{\bar\eta,\mu}$ is a uniform
$\aleph_1$-complete filter on $\mu$ 

\hskip25pt (recalling cf$(\mu) = \theta >
\aleph_0$ as $\mu \in S$) and $\emptyset \notin {\cD}_{\bar\eta,\mu}$.
\end{enumerate}
\mn
[Why?  See $\circledast_9$.]

Now by $\boxplus_3$ we have $\bar\eta \in \Xi \wedge 
e \in \bbB^c_\lambda \cap \bbD_* \Rightarrow Y^\mu_{\bar\eta,e} 
\in {\cD}_{\bar\eta,\mu}$ so
recalling $\circledast_8$ we have $e \in \bbB^c_\lambda \backslash
\bbD_* \Rightarrow Y^\mu_{\bar\eta,e} = \emptyset$ mod 
${\cD}_{\bar\eta,\mu}$ hence
\mn
\begin{enumerate}
\item[$\boxplus_6$]   ${\cP}_{\bar\eta,\mu}
\subseteq \{X \subseteq \mu:X \in {\cD}_{\bar\eta,\mu}$ or
$\mu \backslash X \in {\cD}_{\bar\eta,\mu}\}$.
\end{enumerate}
\mn
Now
\mn
\begin{enumerate}
\item[$\odot_1$]  if $\mu \in S$ then we can find $\bar A^\xi_\mu$ for
$\xi < 2^{2^\mu}$ such that:
\begin{enumerate}
\item[$(a)$]  $\bar A^\xi_\mu = \langle A^\xi_\gamma:\gamma \in
[\mu,2^\mu)\rangle$
\sn
\item[$(b)$]   $A^\xi_\gamma$ is an unbounded subset of $A^*_\mu$
\sn
\item[$(c)$]  $\cD^\xi,\cI^\xi_\mu$ are well defined, i.e. $\emptyset
\in \cD^\xi_\mu$ when we let
\sn
\item[${{}}$]  $(\alpha) \quad {\cD}^\xi_\mu$ be the $\aleph_1$-complete
filter of subsets of $\mu$ generated by 

\hskip35pt  $\{A^\xi_\gamma \backslash
\beta:\gamma \in [\mu,2^\mu)$ and $\beta < \mu\}$

\hskip35pt so $A^*_\mu \backslash \beta \in \cD^\xi_\mu$ for $\beta < \mu$
\sn
\item[${{}}$]  $(\beta) \quad {\cI}^\xi_\mu = \{\mu \backslash B:
B \in {\cD}^\xi_\mu\}$, i.e. the dual ideal
\sn
\item[$(d)$]  moreover if $\xi^1 \ne \xi^2$ are $< 2^{2^\mu}$, then
\[
\{A^*_\mu \backslash A^{\xi^1}_\gamma:\gamma \in [\mu,2^\mu)\} \nsubseteq 
{\cD}^{\xi^2}_\mu \cup {\cI}^{\xi^2}_\mu.
\]
\end{enumerate}
\end{enumerate}
\mn
[Why $\odot_1$ holds?  As $|A^*_\mu| = |\mu|$ is a strong limit
cardinal of cofinality $\theta > \aleph_0$ clearly $\mu = |A^*_\mu| =
|A^*_\mu|^{\aleph_0}$ hence by \cite{EK} there is a sequence
$\langle B_\gamma:\gamma \in [\mu,2^\mu)\rangle$ of subsets
of $A^*_\mu$ such that any  non-trivial Boolean combination of
countably many of them has cardinality $\mu$.  Let $\langle
U_\xi:\xi< 2^{2^\mu}\rangle$ be a sequence of pairwise distinct
subsets of $[\mu,2^\mu)$ each of cardinality $2^{|\mu|}$ no one
included in another and let $\langle 
A^\xi_\gamma:\gamma \in [\mu,2^{|\mu|})\rangle$
list $\{B_\gamma:\gamma \in U_\xi\}$.

\noindent
Now check.]
\mn
\begin{enumerate}
\item[$\odot_2$]  in $\odot_1$ it follows that
\begin{enumerate}
\item[$(e)$]  for every ${\cP} \subseteq {\cP}(\mu)$
for at most one $\xi < 2^{2^\mu}$ we have
\[
\{A^*_\mu \backslash 
A^\xi_\gamma:\gamma \in [\mu,2^\mu)\} \subseteq {\cP} \subseteq
{\cD}^\xi_\mu \cup {\cI}^\xi_\mu.
\]
\sn
\end{enumerate}
\item[$\odot_3$]  for every $\bar\xi = \langle \xi(\mu):\mu \in
S\rangle \in \Pi\{2^{2^\mu}:\mu \in S\}$ there is $\bar\eta =
\bar\eta_{\bar\xi}$ such that:
\begin{enumerate}
\item[$(a)$]  $\bar\eta_{\bar\xi} \in \Xi$ so $\bar\eta_{\bar\xi} =
\langle \eta_{\bar\xi,\varepsilon}:\varepsilon < \lambda\rangle$
\sn
\item[$(b)$]  if $\mu \in S,\gamma \in [\mu,2^\mu)$ then
$\{\varepsilon \in \lambda:\eta_{\bar\xi,\varepsilon}(\gamma) = 1\} =
A^*_\mu \backslash A^{\xi(\mu)}_\gamma$.
\end{enumerate}
\end{enumerate}
\mn
[Why?  Just read the definition of $\Xi$ in $\circledast_9$ and $\bar
A^\xi_\mu$ in $\odot_1$.]
\mn
\begin{enumerate}
\item[$\odot_4$]  if $\mu \in S$ then $\cD_{\bar\eta_{\bar\xi},\mu}
\cup \cI_{\bar\eta_{\bar\xi},\mu} = \cD^{\xi(\mu)}_\mu \cup
\cI^{\xi(\mu)}_\mu$.
\end{enumerate}
\mn
[Why?  Easy, recalling $\boxplus_5(a)$.]
\mn
\begin{enumerate}
\item[$\odot_5$]   if $\gamma(*) < \lambda^+$ and 
$\bar{\bold P}^\gamma 
= \langle \bold P^\gamma_\mu:\mu \in S\rangle$, for $\gamma
< \gamma(*)$ where $\bold P^\gamma_\mu \subseteq {\cP}({\cP}(\mu))$ has 
cardinality $\le 2^\mu$ for $\mu \in S,\gamma < \gamma(*)$
\then \, we can find $\bar \xi = \langle \xi(\mu):\mu \in S\rangle \in
\Pi\{2^{2^\mu}:\mu \in S\}$ 
such that for every $\gamma < \gamma(*)$ the following set 
is not stationary:  $S_{\bar\eta,\gamma} = \{\mu \in S$: 
for some ${\cP} \in \bold P^\gamma_\mu$ we have
$\{A^{\xi(\mu)}_\gamma:\gamma \in [\mu,2^\mu)\} \subseteq
{\cP} \subseteq {\cD}_{\bar\eta,\mu} \cup {\cI}_{\bar\eta,\mu}\}$.
\end{enumerate}
\mn
[Why? Let $\langle u_\alpha:\alpha < \lambda\rangle$ be an increasing
continuous sequence of subsets of $\gamma(*)$ with union $\gamma(*)$
such that $|u_\alpha| \le |\alpha|$ for $\alpha < \lambda$.  Now
for each $\mu \in S$, the family $\cup\{\bold P^\gamma_\mu:\gamma \in
u_\mu\}$ is a family of $\le |u_\mu| \times 2^\mu$ 
subsets of ${\cP}(\mu)$.

Now by clause (e) of $\odot_1$ 
for each $\mu \in S,\gamma \in u_\mu,{\cP} \in \bold P^\gamma_\mu$ 
let $\xi_{\mu,\gamma,{\cP}} < 2^{2^\mu}$ be such
that: if for some $\xi < 2^{2^\mu}$ we have
$\{A^\xi_\gamma:\gamma \in [\mu,2^\mu)\} \subseteq
{\cP} \subseteq {\cD}^\xi_\mu \cup {\cI}^\xi_\mu$ 
then $\xi_{\mu,\gamma,\cP}$ is the first such $\xi$.  
Choose $\xi(\mu) < 2^{2^\mu}$ which does not belong to 
$\{\xi_{\mu,\gamma,{\cP}}:
\gamma \in u_\mu$ and ${\cP} \in \bold P^\gamma_\mu\}$.

So let $\bar\eta = \bar\eta_{\langle \xi(\mu):\mu \in S\rangle} 
\in \Xi$ be as in $\odot_3$ now $\bar\eta$ is as required by 
$\odot_2, \odot_3, \odot_4$.

\noindent
Let us elaborate, why is it as required in $\odot_5$?

First, clearly $\eta_\varepsilon \in {}^\lambda 2$ for $\varepsilon <
\lambda$.  Second, fix $\gamma < \gamma(*)$, then there is $\alpha <
\lambda$ such that $\gamma \in u_\alpha$, so it suffices to show that,
for any $\mu \in S \backslash \alpha$, we have $\mu \notin
S_{\bar\eta,\gamma}$.  So assume $\cP \in \bold P^\gamma_\mu$
satisfies clause (e) of $\odot_1$, and we should prove that
$\neg[\{A^{\xi(\mu)}_\gamma:\gamma \in [\mu,2^\mu)\} \subseteq \cP
\subseteq \cD_{\bar\eta,\mu} \cup \cI_{\bar\eta,\mu}]$; but if for
some $\xi < 2^{2^\mu}$ we have $\{A^\xi_\gamma:\gamma \in
[\mu,2^\mu)\} \subseteq \cP \subseteq \cD_{\bar\eta,\mu} \cup
\cI_{\bar\eta,\mu}$ then necessarily $\xi = \xi_{\mu,\gamma,\cP} \ne
\xi(\mu)$, contradiction to $\odot_1$.]
\mn
\begin{enumerate}
\item[$\odot_6$]    if $\langle M_\gamma:\gamma \le
\gamma(*)\rangle$ is a $\prec$-increasing continuous and $M_\gamma \in
\text{ EC}_\lambda(T)$ and $\bar b_\alpha \in {}^{\ell g(\bar
y)}(M_0)$ for $\alpha < \lambda$ are such that $\langle
\varphi(x,\bar b_\alpha):\alpha < \lambda\rangle$ is
independent, \then \, we can find $N$ such that
\begin{enumerate}
\item[$(\alpha)$]  $M_{\gamma(*)} \prec N \in \text{ EC}_\lambda(T)$
\sn
\item[$(\beta)$]   if $N \prec N' \in \text{ EC}_\lambda(T)$ and
$\gamma < \gamma(*)$ \then \, \footnote{really any pregiven set of
$\le \lambda$ ``forbidden" inv$^\varphi_6$ is O.K. and can make it
work for inv$^\varphi_6(N_\gamma,N')$ for every $\gamma < \gamma(*)$.}
inv$^\varphi_6(M_\gamma,N') \notin
\{\text{inv}^\varphi_6(M_{\gamma_1},M_{\gamma_2}):\gamma \le
\gamma(*)$ and $\gamma_1 < \gamma_2 \le \gamma(*)\}$.
\end{enumerate}
\end{enumerate}
\mn
[Why?  Without loss of generality the universe of $M_{\gamma(*)}$ is
${\cU}_1 \in [\lambda]^\lambda$ such that $\lambda \backslash {\cU}_1$
has cardinality $\lambda$.  Let $\langle u_\alpha:\alpha <
\lambda\rangle$ be as in the proof of $\odot_5$.

For $\gamma(1) < \gamma(2) \le \gamma(*)$ let $\bold
P^{\gamma(1),\gamma(2)}_\delta = \text{ inv}^5_\varphi
(\delta,\text{id}_{N_{\gamma(2)}},M_{\gamma(1)},M_{\gamma(2)})$,
see Definition \ref{inv.10},
clearly inv$^\varphi_6(M_{\gamma(1)},M_{\gamma(2)}) = \langle 
\bold P^{\gamma(1),\gamma(2)}_\delta:\delta <
\lambda\rangle/{\cD}_\lambda$.  So it is enough\footnote{can demand
$\alpha < \lambda \Rightarrow {}^{\omega >}(\omega(\alpha +1))$ if $(N
\backslash M_{\gamma(*)}) \cap [\omega \alpha,\omega \alpha +
\omega)$ is infinite for every $\alpha < \lambda$.}
to find $N$ and sequence $\langle
a_\alpha:\alpha < \lambda\rangle$ of elements of $N$ such that
$M_{\gamma(*)} \prec N,|N| = \lambda$ and for each $\gamma(0) \le
\gamma(*)$, for every $\mu \in S$ except non-stationarily many, the family

\[
\{\{\gamma < \mu:N \models \varphi[a_\gamma,\bar b]\}:\bar b \in
{}^{\ell g(\bar y)}(M_{\gamma(0)})\}
\]

\mn
is not in $\bold P_\mu := \cup\{\bold P^{\gamma(1),\gamma(2)}_\mu:
\gamma(1) < \gamma(2) \le \gamma(*)$ are from $u_\mu\}$.

We choose $\bar\xi = \langle \xi(\mu):\mu \in S\rangle$
as in $\odot_5$; let $\bar\eta = \bar\eta_{\bar\xi}$, see $\odot_3$,
so recalling $\circledast_3$ clearly 
$\langle \Bbb D_{\eta_\varepsilon}:
\varepsilon < \lambda\rangle$ is well defined.  Now
for each $\varepsilon < \alpha$ letting $F$ be from \ref{2d.7} for the
model $M_{\gamma(*)}$ and the sequence $\langle \varphi(x,\bar
b_\alpha):\alpha < \lambda\rangle$, let $p_\varepsilon \in 
\bold S_\varphi(M_{\gamma(*)})$ be such that
for every $\bar b \in {}^{\ell g(\bar y)}(M_{\gamma(*)})$ we have
$\varphi(x,\bar b) \in p_\varepsilon \Leftrightarrow F(\bar b) \in
\Bbb D_{\eta_\varepsilon}$
so $\neg \varphi(x,\bar b) \in p_\varepsilon \Leftrightarrow F(\bar
b) \notin \Bbb D_{\eta_\varepsilon}$.

So by \ref{2d.14} we can find an 
elementary extension $N$ of $M_{\gamma(*)}$ and
$a_\alpha \in N$ for $\alpha < \lambda$ such that
$a_\alpha$ realizes $p_\alpha$, and \wilog \, $N$ has
universe $\subseteq \lambda$ such that $\lambda \backslash |N|$ has
cardinality $\lambda$.  Concerning inv$^6_\varphi$ our demand concerns
what occurs for a club of $\delta < \lambda$ for this.
Let $E \subseteq E_*$ be a club of $\lambda$ such that $\gamma <
\delta \in E \Rightarrow a_\gamma \in N \cap \delta$.  Now in
$\odot_6(\beta)$ we promise something (given $N$) on ``every $N'$ such
that ...", so let $N \prec N' \in \text{ EC}_\lambda(T)$, and 
\wilog \, the universe of $N'$
is $\subseteq \lambda$ and let $\delta \in S \cap E$.  
For any $\gamma \le \gamma(*)$ by $\odot_5$,i.e. by the choice of
$\bar\xi,\bar\eta_\xi$ above there is a club $E_\gamma \subseteq E$ of
$\lambda$ such that for any $\mu \in S \cap E_\gamma$, the
set $S_{\bar\eta,\mu}$ from $\odot_5$ is disjoint to $E_\gamma$, hence
the set $\cP_{\mu,\gamma}
:= \{\{\gamma < \mu:N' \models \varphi[a_\gamma,\bar b]\}:\bar b \in
{}^{\ell g(\bar y)}(M_\gamma)\}$ does not belong to
$\cup\{\bold P^{\gamma(1),\gamma(n)}_\mu:\gamma(1) < \gamma(2) \le 
\gamma(*)$ are from $u_\mu\}$ so we are done.]
\mn
\begin{enumerate}
\item[$\odot_7$]   if $\langle M_\alpha:\alpha < \lambda^+\rangle$
is as in $\boxtimes'$ then for some club $E$ of $\lambda^+$, we have:
if $\alpha_1 < \alpha_2,\beta_1 < \beta_2$ are from $E$ and $\alpha_2
\ne \beta_2$ then

\[
(M_{\alpha_2},M_{\alpha_1}) \ncong (M_{\beta_2},M_{\beta_1}).
\]
\end{enumerate}

\mn
[Why?  For every $\beta < \lambda^+$ we apply $\odot_6$ to $\langle
M_\alpha:\alpha \le \beta\rangle$ and get $N_\beta$ as there so
$M_\beta \prec N_\beta \in \text{ EC}_\lambda(T)$.  As $M =
\cup\{M_\gamma:\gamma < \lambda^+\}$ is saturated, without loss of
generality $N_\beta \prec M$ hence for some $\xi_\beta < \lambda^+$
we have $N_\beta \prec M_{\xi_\beta}$.

Let $E = \{\delta < \lambda^+:\delta$ a limit ordinal such that $\beta
< \delta \Rightarrow \xi_\beta < \delta\}$.

Let $\alpha_1 < \alpha_2,\beta_1 < \beta_2$ be from $E$ such that
$\alpha_2 \ne \beta_2$ and we shall prove that $(M_{\alpha_2},M_{\alpha_1})$
is not isomorphic to $(M_{\beta_2},M_{\beta_1})$.  By symmetry \wilog
\, $\alpha_2 < \beta_2$ and let $\gamma(*) = \text{ max}
\{\alpha_2,\beta_1\}$ so $\gamma(*) < \beta_2$.  Now we apply
$\odot_6$ with $\langle M_\gamma:\gamma \le
\gamma(*)\rangle,N,N',\gamma,\beta_1,\alpha_2,\alpha_2$ here standing
for $\langle M_\gamma:\gamma \le \gamma(*)\rangle,N,N',\gamma(0),
\gamma(1),\gamma(2)$ there so we are clearly done.  
\end{PROOF}
\bigskip\bigskip

%\bibliographystyle{plain}
%\bibliography{lista,listb,listx,listf,liste}

\end{document}